\documentclass{amsart}

\usepackage{amsfonts}
\usepackage{amssymb}
\newcommand{\N}{{\mathbb N}}

\newcommand{\G}{{\hat G}}

\newcommand{\R}{{\mathbb R}}

\newtheorem{thm}{Theorem}[section]
\newtheorem{lemma}[thm]{Lemma}
\newtheorem{cor}[thm]{Corollary}

\newtheorem*{prop}{Proposition}

\theoremstyle{definition}
\newtheorem{rem}[thm]{Remark}
\begin{document}

\title{SEPARABLE SUBSETS OF GFERF NEGATIVELY CURVED GROUPS}

\author{Ashot Minasyan}
%\date{}

\address{Universit\'{e} de Gen\`{e}ve,
Section de Math\'{e}matiques,
2-4 rue du Li\`{e}vre,
Case postale 64,
1211 Gen\`{e}ve 4, Switzerland}

\email{aminasyan@gmail.com}

\begin{abstract} A word hyperbolic group $G$ is called GFERF if every quasiconvex subgroup coincides with the intersection
of finite index subgroups containing it. We show that in any such group, the product of finitely many quasiconvex subgroups
is closed in the profinite topology on $G$.
\end{abstract}

\thanks{This work was partially supported by the NSF grant DMS \#0245600 of A. Ol'shanskii and M. Sapir.}
\keywords{Word Hyperbolic Groups, Profinite Topology, GFERF}
%2000 MSC:
\subjclass[2000]{Primary 20F67, Secondary 20E26.}

\maketitle

\section{Introduction}
Let $G$ be a finitely generated group. The {\it profinite topology} $\mathcal{PT}(G)$  on $G$ is defined by proclaiming
all finite index
normal subgroups to be the basis of open neighborhoods of the identity element. It is easy to see that $G$
equipped with this topology becomes a topological group. This topology is Hausdorff if and only if $G$ is
residually finite.

A subset $P \subseteq G$ will be called {\it separable} if it is closed in the profinite topology on $G$.
Thus, a subgroup $H \le G$ is separable whenever it is an intersection of finite index subgroups.
The group $G$ is said to be {\it locally extended residually finite} (LERF) if every finitely generated subgroup
$H \le G$ is separable.

A famous theorem of M. Hall states that free groups are LERF.
Among other well-known examples of LERF groups are surface groups and fundamental groups of compact Seifert fibred
$3$-manifolds \cite{Scott}. In \cite{Schupp-Cox_sep} P. Schupp provided certain sufficient conditions for a
Coxeter group to be LERF.
R. Gitik  \cite{Gitik-LERFMan} constructed an infinite
family of LERF hyperbolic groups that are fundamental groups of hyperbolic $3$-manifolds.

In 1991 Pin and Reutenauer \cite{Pin-Reut} conjectured that a product of finitely many finitely generated subgroups
in a free group is separable and listed some possible applications to groups and semigroups.
In 1993 Ribes and Zalesski\v{\i} \cite{Rib-Zal} showed that the statement of this
conjecture is true. Later a similar question was studied in other LERF groups by Coulbois \cite{Coulbois},
Gitik \cite{Gitik1}, Niblo \cite{Niblo1}, Steinberg \cite{Steinberg}  and others.

In particular,  Gitik  in \cite[Thm. 1]{Gitik1} proved that in a LERF hyperbolic group, a product of two quasiconvex
subgroups, one of which is malnormal, is separable.

However, many word hyperbolic groups are not LERF. For example, an ascending HNN-extension of a finite rank free group
is never LERF but very often hyperbolic (see \cite{Kap-MTEFG}). So, it makes sense to use the weaker notion below.

We will say that a (word) hyperbolic group $G$ is GFERF if every quasiconvex subgroup $H\le G$
is separable. The definition of a  GFERF Kleinian group $\Gamma$
was  given by Long and Reid in \cite{L-R_1}: $\Gamma$ is called  {\it geometrically finite extended residually finite}
(GFERF) if each geometrically finite subgroup
$H \le \Gamma$ is separable. Our definition is in the same spirit because in any word hyperbolic group
(more generally, in any automatic group) a subgroup is geometrically finite if and only if it is quasiconvex
(see \cite{Swarup_GF-R}).

Long, Reid  and Agol gave several examples of GFERF groups \cite{L-R_1}, \cite{L-R_2}, \cite{L-R_3}.
Hsu and Wise \cite{Hsu-Wise} proved that certain right-angled Artin groups are GFERF.
Some negatively curved (i.e., word hyperbolic) groups with this property  were studied
by Gitik in \cite{Gitik-LERFMan}. In the paper \cite{Wise-fig_8} Wise provided another large family of GFERF
hyperbolic groups; he also showed that the Figure 8 knot group is GFERF. The fact that this group is LERF
follows from the recent proofs by Agol \cite{Agol-tame} and Calegari-Gabai \cite{Cal-Gab-tame}  of
Marden's "tameness" conjecture. This conjecture provides a new way for obtaining LERF and GFERF groups as fundamental
groups of 3-manifolds.

The main goal of this paper is to prove the following

\begin{thm} \label{prof_hyp} Assume $G$ is a GFERF word hyperbolic group, $G_1,G_2,\dots, G_s$ are quasiconvex
subgroups, $s \in \N$. Then the product $G_1G_2 \cdots G_s$ is separable in $G$.
\end{thm}

Since a finitely generated subgroup of a finite rank free group is quasiconvex, the above theorem
generalizes the result of Ribes and Zalesski\v{\i} \cite{Rib-Zal} and provides an alternative proof of the
conjecture \cite{Pin-Reut}. An application of Theorem \ref{prof_hyp} to the case when $s=2$ and $G_2$ is malnormal
gives the statement of Gitik's theorem \cite[Thm. 1]{Gitik1}.

Our proof of Theorem \ref{prof_hyp} uses geometry of quasigeodesics in negatively curved spaces and basic
properties of quasiconvex subgroups.

A subgroup $H$ of a group $G$ will be called {\it almost malnormal} if for every $x \in G \backslash H$ the intersection
$H \cap xHx^{-1}$ is finite. $H$ is said to be {\it elementary} if it is virtually cyclic. It is well known that
in a hyperbolic group $G$ any element of infinite order belongs to a unique {\it maximal elementary subgroup}. Thus,
any maximal elementary subgroup of $G$ is almost malnormal.

A famous open problem in Geometric Group Theory addresses the existence of a (word) hyperbolic group that is
not residually finite. The author would like
to emphasize the importance of studying GFERF hyperbolic groups through the proposition below
(after completion of this article the author found out that a similar statement has been already proved
by D. Wise in \cite[Rem. 11.8]{Wise-res_fin_polyg}).

\begin{prop}
The following are equivalent.
\begin{itemize}
 \item[1)] There exists a non-residually finite hyperbolic group.

 \item[2)] There is a hyperbolic group $G$ having an almost malnormal quasiconvex subgroup $H$ which is not
     separable.

\end{itemize}
\end{prop}

%\vspace{.3cm}
\begin{proof} Assume the first condition holds. In this case Kapovich and Wise \cite[Thm. 1.2]{Kap-Wise},
and, independently, Ol'shanskii \cite[Thm. 2]{Olsh-Bass-Lub}, proved that there exists a non-trivial hyperbolic
group $G$ which has no proper subgroups of finite index at all. Choose an arbitrary maximal elementary subgroup
$H$ of $G$. Obviously $H$ satisfies the condition 2).

Now, suppose 2) holds. Then, according to a theorem of Kharlampovich and Myasnikov \cite[Thm. 2]{Khar-Myas}, the
double $D=G *_H G$ is a hyperbolic group. If the group $D$ were residually finite then we could apply the theorem
of Long and Niblo \cite[\S 2, Lemma]{Long-Niblo} (see also \cite{Niblo1}) stating that $H$ is separable
in $G$. The latter contradicts our assumptions. Hence, $D$ is not residually finite.
\end{proof}

%\vspace{.3cm}
Presently, the author doesn't know of any examples of hyperbolic groups that are not GFERF. So,
it seems reasonable to ask

\vspace{.3cm}
\noindent {\bf Question.} Does there exist a non-GFERF word hyperbolic group ?

\vspace{.3cm}
As one can see from the Proposition, this question may be quite difficult.

Finally, we note that in the case when a hyperbolic group $G$ is GFERF,
Theorem~\ref{prof_hyp} provides a positive solution for Problem 3.11 posed by D. Wise in \cite{Wise-Rips_constr}.
This problem asks whether the double coset $HK$ is separable if $G$ is residually finite and $H,K \le G$ are separable
quasiconvex subgroups.

\vspace{.3cm}
{\noindent \large \bf Acknowledgements}

The author would like to thank Professors Alexander  Ol'shanskii and Mark Sapir for useful
discussions,  the referee and Dr. Dani Wise for their comments.

\section{Auxiliary information}

Suppose $G$ is a group with a fixed finite symmetrized generating set $\mathcal{A}$.
 If $g\in G$, $|g|_G$ will denote the length of a shortest word over $\mathcal{A}$ representing
$g$. Now we can define the standard left-invariant word metric $d(\cdot,\cdot)$ on $G$ by setting
$d(x,y) \stackrel{def}{=} |x^{-1}y|_G$ for arbitrary $x,y \in G$. This metric extends to a metric on the Cayley
graph $\Gamma(G,\mathcal{A})$ of the group $G$
after endowing every edge with the metric of the segment $[0,1] \subset \R$.

A subset $Q$ of $G$ is said to be $\varepsilon$-quasiconvex (where $\varepsilon \ge 0$) if for any pair of elements
$u,v \in Q$ and any geodesic segment $p$ connecting $u$ and $v$, $p$ belongs to a closed $\varepsilon$-neighborhood
of $Q$ in $\Gamma(G,\mathcal{A})$. A subset $Q \subset G$ is quasiconvex if it is $\varepsilon$-quasiconvex for some
$\varepsilon \ge 0$.

For any two points $x,y \in \Gamma(G,\mathcal{A})$ we fix a geodesic path between them and denote it by $[x,y]$.
If $x,y,w \in \Gamma(G,\mathcal{A})$, then the number
$$(x|y)_w \stackrel{def}{=} \frac12 \Bigl(d(x,w)+d(y,w)-d(x,y) \Bigr)$$
is called the {\it Gromov product} of $x$ and $y$ with respect to $w$.

\begin{rem} \label{prodinvar} Since the metric is left-invariant, for arbitrary $x,y,w \in G$ we have
$(x|y)_w=(w^{-1}x|w^{-1}y)_{1_G}$.
\end{rem}

Let $abc$ be a geodesic triangle in $\Gamma(G,\mathcal{A})$.
There exist "special" points $O_a \in [b,c]$, $O_b \in [a,c]$, $O_c \in [a,b]$ with the properties:
$d(a,O_b) = d(a,O_c) = \alpha$, $d(b,O_a) = $ $=d(b,O_c) = \beta$, $d(c,O_a) = d(c,O_b) = \gamma$. It is easy to see
that $\alpha = (b|c)_a$, $\beta = (a|c)_b$, $\gamma = (a|b)_c$. Two points
$O \in [a,b]$ and $O' \in [a,c]$ are called $a$-{\it equidistant} if $d(a,O) = d(a,O') \le \alpha$.
The triangle $abc$ is said to be $\delta$-{\it thin} if for any two points $O,O'$ lying on its sides and
equidistant from one of its vertices, $d(O,O') \le \delta$ holds.

The group $G$ is said to be (word) hyperbolic (or negatively curved) if there is $\delta\ge 0$ such that
every geodesic triangle in $\Gamma(G,\mathcal{A})$ is $\delta$-thin (for more theory the reader is referred to \cite{Ghys},\cite{Mihalik}).

For a hyperbolic group $G$, the property of a subset to be quasiconvex
does not depend on the choice of a generating set $\mathcal{A}$ (see \cite{Gromov}). A quasiconvex subgroup
of a finitely generated group is finitely generated itself (\cite{Mihalik}, \cite{Short}). A conjugate of a quasiconvex
subgroup is quasiconvex as well (\cite[Remark 2.2]{hyp}).

%\begin{lemma} \label{qc-fg,inter}
%{\normalfont (\cite{Mihalik}, \cite{Short})} Suppose $A$ is a quasiconvex subgroup
%of a group $G$ with a fixed finite generating set $\mathcal{A}$.
%Then $A$ is finitely generated. %If $B \le G$ is another quasiconvex subgroup, then $A \cap B$ is also quasiconvex
%\end{lemma}

Fix an arbitrary GFERF hyperbolic group $G$. Then for $n \in \N$, $f_0,f_1,$ $\dots,$ $f_n \in G$ and any
quasiconvex subgroups
$G_1,\dots,G_n \le G$, the subset
\begin{equation} \label{q-c_prod} P =f_0G_1f_1G_2 \cdot \dots \cdot f_{n-1}G_nf_n
\end{equation} is called a {\it quasiconvex product} (here we use the terminology from \cite{hyp}). Such a subset is always
quasiconvex (\cite[Prop. 3.14]{Zeph}, \cite[Cor. 2.1]{hyp}).
%
%We can define a collection of subsets ${\mathcal{K}}(G) \subset 2^G$ to be the minimal collection containing the empty set,
%all quasiconvex products and closed under finite unions.
%
%Then for any pair of subsets $U_1,U_2 \in {\mathcal{K}}(G)$ we have $U_1\cdot U_2 \in \mathcal{K}(G)$ and
%$U_1 \cap U_2 \in \mathcal{K}(G)$ (see \cite{hyp}).

\begin{rem} \label{transform} Assume that $n \in \N$ and for any $n$ quasiconvex subgroups of the group $G$,
their product is closed in $\mathcal{PT}(G)$. Then any quasiconvex product $P$ defined by
(\ref{q-c_prod}) is also closed in $\mathcal{PT}(G)$.
\end{rem}

Indeed, %since $G$ endowed with $\mathcal{PT}(G)$ is a topological group, it is enough to
observe that $P=f\G_1\cdot \dots \cdot \G_n$ where $f=f_0f_1\cdot \dots \cdot f_n \in G$
and $\G_i=(f_if_{i+1} \cdots f_n)^{-1}G_i(f_if_{i+1} \cdots f_n)$ -- quasiconvex subgroups of $G$.
By the assumptions, $\G_1\cdot \dots \cdot \G_n$ is separable, and since $G$ (endowed with $\mathcal{PT}(G)$)
is a topological group, left translation by the element $f^{-1} \in G$ is a continuous operation, hence
$P$ is also separable.

\begin{lemma} \label{smallgromprod} Assume that $G$ is a $\delta$-hyperbolic group with respect to a
finite generating set $\mathcal{A}$  and $A$, $B$ are $\varepsilon$-quasiconvex subgroups.
There exists a constant $C_0=C_0(\delta,\varepsilon,G,\mathcal{A}) \ge 0$ such that for any $a \in A$, $b \in B$
the inequality
$(a^{-1}|b)_{1_G} \le C_0$  holds whenever $a$ is a shortest representative of the coset $a(A\cap B)$.
\end{lemma}

\begin{proof} Define a finite subset of the group $G$ by $\Theta=\{g \in AB~|~|g|_G \le 2\varepsilon+\delta\}$.
For every $g \in \Theta$ choose a pair $(x,y)\in A\times B$ satisfying $g =x^{-1}y$; let $\Omega \subset A\times B$
denote the (finite) set of these pairs. Consider
$$\Omega_1=\{x \in A~|~(x,y) \in \Omega \mbox{ for some } y \in B\}.$$

Then one can define the number $C_0=\max\{|x|_G~|~x \in \Omega_1\}+\varepsilon <\infty$.

Now, assume that $(a^{-1}|b)_{1_G} > C_0$, for some $a \in A$, $b \in B$
where $a$ is a shortest representative of the coset $a(A\cap B)$. Let  $\alpha$ and $\beta$ denote
the "special" points of the triangle $1_Ga^{-1}b$ (in $\Gamma(G,\mathcal{A})$)
on the sides $[1_G,a^{-1}]$ and $[1_G,b]$ respectively.
Since $A$ and $B$ are $\varepsilon$-quasiconvex
there are elements $a_1\in A$ and $b_1 \in B$ that are $\varepsilon$-close to $\alpha$ and $\beta$ correspondingly.
Using the triangle inequality and $\delta$-hyperbolicity of the space $\Gamma(G,\mathcal{A})$ we obtain
\begin{multline*} |aa_1|_G=d(a^{-1},a_1)\le d(a^{-1},\alpha)+\varepsilon = d(a^{-1},1_G)-d(\alpha,1_G)+\varepsilon
= \\ d(a^{-1},1_G)-(a^{-1}|b)_{1_G} +\varepsilon<|a|_G-C_0+\varepsilon,\end{multline*}
$$  |a_1^{-1}b_1|_G=d(a_1,b_1)\le d(a_1,\alpha)+d(\alpha,\beta)+d(\beta,b_1)\le 2\varepsilon+\delta.$$

By definition, there exists a pair of elements $(x,y)\in \Omega$ with $a_1^{-1}b_1=x^{-1}y$, thus
$a_1x^{-1}=b_1y^{-1}\in A\cap B$.
Now, $a(a_1x^{-1}) \in a(A \cap B)$ and this element is shorter than $a$ because
$$|aa_1x^{-1}|_G \le |aa_1|_G+|x|_G < |a|_G-(C_0-\varepsilon-|x|_G) \le |a|_G.$$
Thus we achieve a contradiction with our assumptions.
\end{proof}

Let $p$ be a path in the Cayley graph of $G$. Then $p_-$, $p_{+}$  will denote the initial and
the final points of $p$, $||p||$ -- its length. We will use $elem(p)$ to denote the element of the group $G$ represented by the
word written on $p$. A path $q$ is called $(\lambda,c)$-{ \it quasigeodesic} if there exist
$0<\lambda \le 1$, $c \ge 0$, such that
for any subpath $p$ of $q$ the inequality $\lambda ||p|| - c \le d(p_-,p_+)$ holds.

The statement below is an analog of the fact that in a negatively curved space $k$-local geodesics
are quasigeodesics for any sufficiently large $k$.

\begin{lemma}  \label{quasigeod}  {\normalfont \cite[Lemma 4.2]{paper3}}
Let $\bar \lambda >0$, $\bar c \ge 0$, $C_0 \ge 14 \delta$,
$C_1 = 12(C_0+\delta)+\bar c + 1$ be given.
Then  for $\lambda = \bar \lambda/4>0$ there exist $c=c(\bar \lambda,\bar c,C_0) \ge 0$ satisfying the statement below.

Assume $N \in \N$, $x_i \in \Gamma(G,\mathcal{A})$, $i=0,\dots,N$, and $q_i$ are $(\bar \lambda,\bar c)$-quasigeodesic paths between
$x_{i-1}$ and $x_i$ in $\Gamma(G,\mathcal{A})$, $i=1,\dots,N$. If $\|q_i\| \ge (C_1+ \bar c)/{\bar \lambda}$, $i=1,\dots,N$, and
$(x_{i-1}|x_{i+1})_{x_i} \le C_0$ for all $i=1,\dots,N-1$,
then the path $q$ obtained as a consecutive concatenation of
$q_1,q_2,\dots,q_N$ is $(\lambda,c)$-quasigeodesic.
\end{lemma}

For any element $x \in G$ and $N\ge 0$ the closed ball centered at $x$ of radius $N$ will be denoted by
$\mathcal{O}_N (x) = \{y\in G~|~d(x,y) \le N \}$.

%If $Q \subset G$, $N \ge 0$, the closed $N$-neighborhood of $Q$ will be denoted by
%$$\mathcal{O}_N (Q) \stackrel{def}{=} \{x\in G~|~d(x,Q) \le N \}.$$

\begin{lemma} \label{mainlemmafor2} Assume $G$ is a $\delta$-hyperbolic group, $A$ and $B$ are
$\varepsilon$-quasi\-convex
subgroups. Then for any $N \ge 0$ there exists $N_1=N_1(N,\delta,\varepsilon,G,\mathcal{A})\ge 0$ such that the following holds.
Suppose the subgroups $A' \le A$ and $B' \le B$ satisfy
$A\cap B =A'\cap B'$, $\mathcal{O}_{N_1}(1_G) \cap (A' \cup B') \subset A\cap B$. Then for the subgroup
$H=\langle A',B'\rangle \le G$ one has $$ \mathcal{O}_N(1_G) \cap A H  B \subset AB~. $$
\end{lemma}

\begin{proof} First, let $\hat C_0= \hat C_0(\delta,\varepsilon,G,\mathcal{A})$ be the constant given by Lemma
\ref{smallgromprod}. Define $C_0=\max\{\hat C_0,14\delta\}$, $\bar \lambda=1$, $\bar c=0$ and
$C_1 = 12(C_0+\delta)+\bar c + 1$. Now apply Lemma \ref{quasigeod} to find $\lambda = \bar \lambda/4=1/4>0$ and
$c=c(\bar \lambda,\bar c,C_0) \ge 0$ from its claim.

Set $N_1=(N+c+2C_1)/\lambda$ ~and let $A' \le A$ and $B' \le B$ satisfy the conditions of the lemma. Thus,
%Since the intersection $A\cap B$ is quasiconvex (Lemma \ref{qc-fg,inter}) and
%$G$ is GFERF,
%we can find finite index subgroups $A' \le_f A$ and $B' \le_f B$ satisfying $A \cap B \le A'$, $A \cap B \le B'$ and
\begin{equation} \label{A',B'}
A' \cap \mathcal{O}_{N_1}(1_G)\subset A\cap B,~ B' \cap \mathcal{O}_{N_1}(1_G)\subset A\cap B
\end{equation}
%As it follows the construction, $A'$ and $B'$ satisfy all the claims of the lemma except, possibly, $(*)$.
%
%Now, $\displaystyle A=\bigsqcup_{i=1}^ma_i\A$, $\displaystyle B=\bigsqcup_{j=1}^n\B b_j$ ~for some $m,n \in \N$,
%$a_i \in A$, $b_j \in B$ for all $i,j$. Set
%$\displaystyle N_1=N+\max_{1\le i \le m}|a_i|_G + \max_{1\le j \le n}|b_j|_G$.

%Again, using the GFERF property of $G$, we are able to find finite index subgroups $A'\le_f A$ and
%$B' \le_f B$ with $A \cap B \le A'$, $A \cap B \le B'$ and
%\begin{equation} \label{A',B'}
%A' \cap \mathcal{O}_{(N_1+c+C_1)/\lambda}(1_G)\subset A\cap B,~ B' \cap \mathcal{O}_{(N_1+c+C_1)/\lambda}(1_G)\subset A\cap B.
%\end{equation}

Define the subgroup $H=\langle A',B'\rangle \le G$ and consider an arbitrary element \\
$g\in A H B \backslash (AB)$.

%{\bf Case 1.} Suppose $g \in A H B$.
Then \begin{equation} \label{expr_for_g} g=x_0y_1x_1y_2\cdots x_{l}y_{l+1},\end{equation}
where $l \in \N \cup\{0\}$,
$x_0 \in A$, $x_i \in A'\backslash \{1_G\}$, $y_i \in B'\backslash \{1_G\}$, $i=1,\dots,l$, $y_{l+1} \in B$.
Moreover, we can assume that $x_0,x_1,\dots,x_l,y_1,\dots,y_l$ are shortest representatives of their left
cosets modulo $A \cap B$ (indeed, if there is $\tilde x_0 = x_0z$ with $z \in A\cap B$ and $|\tilde x_0|_G<|x_0|_G$,
then $\tilde x_0 \in A$, $g=\tilde x_0 (zy_1)x_1y_2\cdots x_ly_{l+1}$ where $zy_1 \in B'$
because of the construction of $A'$; and then a similar procedure can be performed for $zy_1$, and so on) and $l$
is the smallest such integer. Therefore
\begin{equation} \label{x_i,y_i} x_i \in A'\backslash (A\cap B),~y_i \in B'\backslash (A\cap B),~i=1,\dots,l.
\end{equation}

Observe that since $g \notin AB$, $l \ge 1$ and $y_1 \in B'\backslash (A\cap B)$.
Choose geodesic paths $q_1$, $q_2$, $\dots$, $q_{2l+2}$ in
$\Gamma(G,\mathcal{A})$ as follows: $(q_1)_-=1_G$, $elem(q_1)=x_0$, $(q_2)_-=(q_1)_+$, $elem(q_2)=y_1$,
$\dots$,
%$(q_{2l})_-=(q_{2l-1})_+$, $elem(q_{2l})=y_l$; now, if $y_{l+1} \in A\cap B$
%then $x_ly_{l+1} \in A'$
$(q_{2l+2})_-=(q_{2l+1})_+$, $elem(q_{2l+2})=y_{l+1}$. Thus, $(q_{2l+2})_+=g$. Using (\ref{x_i,y_i}) and
(\ref{A',B'}) we obtain
$\|q_i\| > N_1 \ge C_1=(C_1+ \bar c)/{\bar \lambda}$, $i=2,3,\dots,2l+1$, and $((q_i)_-|(q_{i+1})_+)_{(q_i)_+} \le C_0$
(by Remark \ref{prodinvar} and Lemma \ref{smallgromprod}) for $i=1,\dots,2l+1$.

Now, there can occur four different situations depending on how long the paths $q_1$ and $q_{2l+2}$ are:
\begin{itemize}
    \item[\bf (a)] $\|q_1\|< C_1$ and $\|q_{2l+2}\|<C_1$;
    \item[\bf (b)] $\|q_1\| \ge C_1$ and $\|q_{2l+2}\|<C_1$;
    \item[\bf (c)] $\|q_1\|< C_1$ and $\|q_{2l+2}\| \ge C_1$;
    \item[\bf (d)] $\|q_1\| \ge C_1$ and $\|q_{2l+2}\| \ge C_1$.
\end{itemize}

Let us consider the situation (b) (the others can be resolved in a completely analogous fashion).
Then the path $q=q_1q_2\dots q_{2l+1}$ satisfies all the conditions of Lemma \ref{quasigeod}, hence
it is $(\lambda,c)$-quasigeodesic (for the numbers $\lambda,c$
defined in the beginning of the proof). Recalling (\ref{A',B'}) we get
$$|g|_G \ge d(q_-,q_+) -d(q_+,g) \ge \lambda \|q\| -c -\|q_{2l+2}\|\ge \lambda\|q_2\| -c -C_1$$
$$=\lambda|y_1|_G -c -C_1> \lambda((N+c+2C_1)/\lambda)-c-C_1\ge N~.$$

Similarly, one can show that $|g|_G>N$ in the other three situations.

Thus, we have $AHB\cap \mathcal{O}_{N}(1_G) \subset AB$ and the lemma is proved. \end{proof}

%Let $G$ be a $\delta$-hyperbolic group, $A$, $B$ -- its $\varepsilon$-quasiconvex subgroups.
%Suppose $A' \le A$, $B' \le B$ satisfy $A\cap B = A' \cap B'$. Then any element
%$h \in A\langle A',B'\rangle B$ can be written as follows:
%\begin{equation} \label{red_def} h=x_0y_1x_1y_2\cdots x_{l}y_{l+1}, \end{equation} where $l \in \N \cup\{0\}$,
%$x_0 \in A$, $x_i \in A'\backslash \{1_G\}$, $y_i \in B'\backslash \{1_G\}$, $i=1,\dots,l$, $y_{l+1} \in B$.
%
%\vspace{.3cm}
%\underline{\bf Definition.} The representation (\ref{red_def}) will be called {\it reduced} if the elements
% $x_0,x_1,\dots,x_l,y_1,\dots,y_l$ are shortest representatives of their left
%cosets modulo $A \cap B$.
%
%\vspace{.3cm}
%As we saw during the proof of Lemma \ref{mainlemmafor2} any element $h \in A\langle A',B'\rangle B$ possesses at least
%one reduced representation. Also, from this proof one immediately achieves
%
%\begin{cor} There exists constants $M \ge 0, \lambda >0$ and $c \ge 0$ (depending on $G$, $\delta$ and $\varepsilon$)
%such that if all elements in $A'$ and $B'$ which are shorter than $M$ belong to $A\cap B$
%then the following property holds.
%
%For an arbitrary reduced representation (\ref{red_def}) consider geodesic paths $q_1,q_2,\dots,$ $q_{2l+2}$ in the
%Cayley graph
%$\Gamma(G,\mathcal{A})$ with $(q_{i+1})_-=(q_{i})_+$, $i=0,\dots,2l$, and $elem(q_1)=x_0$, $elem(q_2)=y_1$, \dots, $elem(q_{2l+1})=x_l$,
%$elem(q_{2l+2})=y_{l+1}$. Then the broken line $q=q_1q_2\dots q_{2l+2}$ is $(\lambda,c)$-quasigeodesic.
%\end{cor}

%\vspace{.3cm}
Note that during the proof of Lemma \ref{mainlemmafor2} for each $g \in H=\langle A',B'\rangle$ we constructed a
presentation (\ref{expr_for_g}) and a corresponding quasigeodesic path $q=q_1 \dots q_{2l+2}$ connecting $1_G$ and
$g$ in $\Gamma(G,\mathcal{A})$. Since geodesics and quasigeodesics with same ends are mutually close (\cite[3.3]{Mihalik}), the geodesic
$[1_G,g]$ will lie in some neighborhood of $q$. If, in addition, the subgroups $A'$ and $B'$ are
$\varepsilon'$-quasiconvex, $q$ will belong to a closed $\varepsilon'$-neighborhood of $H$ in $\Gamma(G,\mathcal{A})$.
Hence $H$ becomes quasiconvex itself. Thus, one obtains the statement below, first proved by R. Gitik:

\begin{lemma} {\normalfont (\cite[Thm. 1]{Gitik-pong})} \label{AB-qc} Let $A$ and $B$ be $\varepsilon$-quasiconvex
subgroups of a $\delta$-hyperbolic group $G$. There exists a constant $C_2$, which depends only on $G$, $\delta$ and
$\varepsilon$, with the following property.
For any quasiconvex subgroups $A' \le A$ and $B' \le B$ with $A' \cap B'=A\cap B$, if all elements in $A'$ and $B'$
shorter than $C_2$ belong to $A \cap B$, then the subgroup $\langle A',B' \rangle$ is also quasiconvex in $G$.
\end{lemma}

\begin{cor} \label{AB-closed} If $G$ is a GFERF hyperbolic group and $A$, $B$ are its quasiconvex subgroups then
the double coset $AB$ is separable in $G$.
\end{cor}

\begin{proof} It is enough to show that for arbitrary $g \in G \backslash (AB)$ there exists a closed
(in the profinite topology) subset $K$ of $G$ such that $AB \subseteq K$ and $g \notin K$.
Let $C_2$ be the constant given by Lemma \ref{AB-qc}.
Set $N=|g|_G$ and find the corresponding $N_1 \ge 0$ from the claim of Lemma \ref{mainlemmafor2}. Denote
$N_2=\max\{N_1,C_2\}$. Since the subgroups $A$ and $B$ are closed in $\mathcal{PT}(G)$, then so is $A\cap B$; therefore
%is quasiconvex (Lemma \ref{qc-fg,inter}) and $G$ is GFERF,
there exist subgroups $A' \le_f A$ and $B' \le_f B$ (having finite indices in $A$ and $B$ correspondingly)
such that $A\cap B \subset A'$, $A\cap B \subset  B'$ and $\mathcal{O}_{N_2}(1_G) \cap (A' \cup B') \subset A\cap B$.
Applying Lemma \ref{mainlemmafor2} to $H=\langle A',B'\rangle \le G$  we achieve $g \notin A H B$.

Now, $\displaystyle A=\bigsqcup_{i=1}^ma_iA'$, $\displaystyle B=\bigsqcup_{j=1}^nB' b_j$ ~for some $m,n \in \N$,
$a_i \in A$, $b_j \in B$ for all $i,j$. Since a finite index subgroup of a quasiconvex subgroup is itself quasiconvex,
$H$ is quasiconvex by Lemma \ref{AB-qc}, hence it is closed in $\mathcal{PT}(G)$ as $G$ is GFERF.
Therefore the sets $a_iHb_j$ are
closed for any $i,j$, and, consequently, their finite union
$$ K\stackrel{def}{=}\bigcup_{i=1}^m\bigcup_{j=1}^n a_iHb_j$$ is
closed too. It remains to observe that $AB \subset K= AHB$, thus $g\notin K$. Q.e.d. \end{proof}

\section{Proof of Theorem \ref{prof_hyp}}
{\renewcommand{\i}{^{(i)}}

\begin{proof} We will use induction on $s$. If $s=1$, the statement follows from the definition of a GFERF group.
The case $s=2$ is given by Corollary \ref{AB-closed}. So, we can now assume that $s>2$ and the statement is already
proved for a product of any $(s-1)$ quasiconvex subgroups.

For our convenience, denote $k=s-2$, $A=G_{s-1}$, $B=G_s$. Let $\{A_i~|~i\in \N\}$, $\{B_i~|~i\in \N\}$ be enumerations
of all finite index subgroups containing $A \cap B$ in $A$ and $B$ correspondingly. Define the sequences
$$A\i = \bigcap_{j=1}^i A_j,~B\i = \bigcap_{j=1}^i B_j~.$$

Now, due to the construction, $A\cap B \subset A\i \le_f A$ and $A\cap B \subset B\i \le_f B$ for all $i$.
And (as we saw in the proof of Corollary \ref{AB-closed}) for every $i \in \N$ there are $m=m(i),n=n(i) \in \N$
and elements $a_1,\dots,a_m \in A$, $b_1,\dots, b_n \in B$ such that
\begin{equation} \label{ahb} A\langle A\i,B\i \rangle B= \bigcup_{p=1}^m\bigcup_{r=1}^n a_p\langle A\i,B\i \rangle b_r~.
\end{equation}

% Also we can make

\begin{rem} \label{seq_of_subgrps}
For any finite index subgroup $H$ of $G$ satisfying $A \cap B \le H$ there exists
$I \in \N$ such that $A\i,B\i \le H$ for all $i \ge I$.
\end{rem}

Since $A$ and $B$ are separable in $G$, their intersection $A \cap B$ is separable as well, and we have
$$A\cap B =\bigcap_{i=1}^\infty A\i=\bigcap_{i=1}^\infty B\i.$$

Without loss of generality, we can assume that the subgroups $G_1,\dots,G_k$,$A$, $B$ are $\varepsilon$-quasiconvex
for a fixed $\varepsilon \ge 0$.
Let $\hat C_0= \hat C_0(\delta,\varepsilon,G,\mathcal{A})$ be the constant given by Lemma
\ref{smallgromprod}. Define $C_0=\max\{\hat C_0,14\delta\}$, $\bar \lambda=1$, $\bar c=0$ and
$C_1 = 12(C_0+\delta)+\bar c + 1$. Now apply Lemma \ref{quasigeod} to find $\lambda = \bar \lambda/4>0$ and
$c=c(\bar \lambda,\bar c,C_0) \ge 0$ from its claim.

Let $C_2=C_2(\delta,\varepsilon,G)$ be the constant from the claim of Lemma \ref{AB-qc}. Since the group $G$ is
GFERF, there exist $A' \le_f A$ and $B' \le_f B$ such that $A'\cap B'=A\cap B$ and all the elements in
$A'$ and $B'$ shorter than $C_2$ belong to $A\cap B$. Therefore,  we can find an
index $I_1 \in \N$ such that $A\i \le A'$, $B\i \le B'$ for all $i\ge I_1$, hence, according to Lemma \ref{AB-qc},
the subgroup $\langle A\i,B\i \rangle \le G$ is quasiconvex.

Arguing by contradiction, suppose there exists  $g \in G \backslash(G_1G_2\dots G_kAB)$
which belongs to the closure of $G_1G_2\dots G_kAB$ in $\mathcal{PT}(G)$. Keeping in mind formula (\ref{ahb}) and Remark
\ref{transform}, for any $i\ge I_1$ we can apply the induction hypothesis to the product
$$P_i\stackrel{def}{=}G_1G_2\cdots G_kA\langle A\i,B\i \rangle B$$ to show that it is closed in $\mathcal{PT}(G)$.

Obviously, $G_1G_2\dots G_kAB \subseteq P_i$, hence $g \in P_i$ for every $i \ge I_1$.
Thus, for each $i\ge I_1$ one can find $l=l(i) \in \N \cup\{0\}$ and elements $z_1\i \in G_1$, $\dots$, $z_k\i \in G_k$,
$x_0\i \in A$, $x_j\i \in A\i \backslash (A\cap B)$, $y_j\i \in B\i \backslash (A\cap B)$, $i=1,\dots,l$,
$y_{l+1}\i \in B$ satisfying
\begin{equation} \label{form_of_g} g=z_1\i \cdot \dots \cdot z_k\i x_0\i y_1\i x_1\i \cdot \dots\cdot x_l\i y_{l+1}\i.
\end{equation}
Moreover, as in the proof of Lemma \ref{mainlemmafor2}, we can assume that $z_t$ is a shortest representative of
its left coset modulo $G_t \cap G_{t+1}$ for $t=1,\dots,k-1$, $z_k$ is a shortest representative of
its left coset modulo $G_k \cap A$, and $x_0,x_j,y_j$ are shortest representatives of
their left cosets modulo $A \cap B$ for $j=1,\dots,l$.

Now we have to consider several possibilities.

{\bf CASE 1.} For some $t \in \{1,\dots ,k\}$ we have $\displaystyle \liminf_{i\to \infty} |z_t\i|_G <\infty$.

Then, by passing to a subsequence, we can assume that $z_t\i=z_t \in G_t$ for all $i$. Using (\ref{form_of_g})
and our assumptions on $g$ we obtain
\begin{equation} \label{case1} g \in G_1\cdots G_{t-1}z_tG_{t+1}\cdots G_k
A\langle A\i,B\i \rangle B, \end{equation}
$$\mbox{and } g \notin G_1\cdots G_{t-1}z_tG_{t+1}\cdots G_kA B ~\mbox{ for all $i$.}$$

By Remark \ref{transform} and the induction hypothesis, the subset
$$G_1\cdots G_{t-1}z_tG_{t+1}\cdots G_kA B$$ is closed in $\mathcal{PT}(G)$, consequently, there exists a
 normal subgroup $K$ of finite index in $G$ such that
 $$gK \cap G_1\cdots G_{t-1}z_tG_{t+1}\cdots G_kA B =\emptyset.$$
Since $BK=KB=KBB$,
\begin{equation} \label{case1-2}g \notin G_1\cdots G_{t-1}z_tG_{t+1}\cdots G_kA K B=
G_1\cdots G_{t-1}z_tG_{t+1}\cdots G_kA H B, \end{equation}
where $H=KB$ is a finite index subgroup of $G$ containing $A\cap B$. Applying Remark \ref{seq_of_subgrps}
we achieve that $\langle A\i,B\i \rangle \le H$ for every sufficiently large $i$, thus
\begin{equation} \label{case1-3}G_{t-1}z_tG_{t+1}\cdots G_kA \langle A\i,B\i \rangle B \subseteq G_{t-1}z_tG_{t+1}\cdots
 G_kA H B. \end{equation}

Combining (\ref{case1}),(\ref{case1-2}) and (\ref{case1-3}) together we obtain a contradiction.

{\bf CASE 2.}
Suppose $\displaystyle \liminf_{i \to \infty} |x_0\i|_G<\infty$.

Again, by passing to a subsequence, we are able to assume that $x_0\i=x_0 \in A$ for all $i$. Thus,
\begin{equation} \label{case2} g \in G_1\cdot \dots\cdot G_kx_0\langle A\i,B\i \rangle B ~\mbox{ for all $i$.}
\end{equation}
 Now, since the subset $G_1\cdots G_kx_0B$ is closed in $\mathcal{PT}(G)$,
we can find a normal subgroup $K$ having finite index in $G$ and satisfying
\begin{equation} \label{case2-2} g \notin G_1\cdots G_kx_0KB=G_1\cdots G_kx_0HB,
\end{equation} where $H=KB \le_f G$ and $A \cap B \le H$. Similarly to Case 1, formula (\ref{case2-2})
leads to a contradiction with formula (\ref{case2}).

{\bf CASE 3.} Suppose $\displaystyle \liminf_{i \to \infty} |y_{l+1}\i|_G<\infty$ (though $l$ may depend on $i$,
it doesn't matter for us).

This case can be resolved in the same way as Case 2.

\noindent And, finally, the last

{\bf CASE 4.} For every $t\in \{1,\dots ,k\}$ we have $\displaystyle \lim_{i\to \infty} |z_t\i|_G =\infty$ and,
in addition, $\displaystyle \lim_{i\to \infty} |x_0\i|_G =\lim_{i\to \infty} |y_{l+1}\i|_G=\infty$.

Then for some $i>I_1$, we will have $|z_t\i|_G > C_3$ for $t=1,\dots,k$, $|x_0\i|_G >C_3$,
$|x_j\i|_G >C_3$ (since $x_j\i \in A\i \backslash (A\cap B)$ and $A\cap B= \bigcap_{i=1}^\infty A\i$),
$|y_j\i|_G >C_3$ for $j=1,\dots,l$, $|y_{l+1}\i|_G >C_3$,
where $$C_3 \stackrel{def}{=} \max\left\{\frac{C_1+\bar c}{\bar \lambda},\frac{|g|_G+ c}{\lambda}\right\}.$$

Choose the geodesic paths $q_1,\dots,q_{k+2l+2}$ in $\Gamma(G,\mathcal{A})$ as follows: $(q_1)_-=1_G$, $elem(q_1)=z_1\i$,
$\dots$, $(q_k)_-=(q_{k-1})_+$, $elem(q_k)=z_k\i$, $(q_{k+1})_-=(q_{k})_+$, \\$elem(q_{k+1})=x_0\i$,
$(q_{k+2})_-=(q_{k+1})_+$, $elem(q_{k+2})=y_1\i$, $\dots$, $(q_{k+2l+2})_-=(q_{k+2l+1})_+$,
$elem(q_{k+2l+2})=y_{l+1}\i$. Recalling (\ref{form_of_g}) we see that $(q_{k+2l+2})_+=g$.

Now, by the construction of presentation (\ref{form_of_g}), we can first apply Lemma \ref{smallgromprod} and then
Lemma \ref{quasigeod} to the broken line $q=q_1\dots q_{k+2l+2}$. Thus, $q$ is $(\lambda,c)$-quasigeodesic.
Since $q_-=1_G$, $q_+=g$, we get
$$|g|_G=d(q_-,q_+) \ge \lambda \|q\|-c \ge \lambda \|q_1\|-c>\lambda C_3-c \ge |g|_G~.$$

The contradiction achieved finishes the proof. \end{proof} }

\end{document}